\documentclass[12pt]{elsarticle}

\usepackage{mystyleElsArticle}
\usepackage{amsthm}
\usepackage{graphicx}
\usepackage{caption}
\usepackage{subcaption}

\newcommand{\rev}[1]{\textcolor{black}{#1}}

\newcommand{\oM}{\mathcal{M}}
\newcommand{\oS}{\mathcal{S}}
\newcommand{\oL}{\mathcal{L}}
\newcommand{\oLt}{\tilde{\oL}}

\newcommand{\oQ}{\mathcal{Q}}
\newcommand{\oQt}{\tilde{\oQ}}

\newcommand{\inp}[2]{\left\langle #1 , \, #2 \right\rangle}

\newcommand{\vU}{\vec{U}}

\newcommand{\vY}{\vec{Y}}

\newcommand{\mLa}{\mat{\Lambda}}
\newcommand{\mC}{\mat{C}}

\newcommand{\mD}{\mat{D}}

\newcommand{\mI}{\mat{I}}
\newcommand{\mA}{\mat{A}}

\newcommand{\mM}{\mat{M}}

\newcommand{\hmpn}{\mathbf{H}}

\newcommand{\mindex}[1]{\boldsymbol{#1}}
\newcommand{\minda}{\mindex{\alpha}}

\newcommand{\gvec}[1]{\boldsymbol{#1}}
\newcommand{\vxi}{\gvec{\xi}}

\newcommand{\mset}[1]{\mathfrak{#1}}

\newcommand{\sJ}{\mset{J}}

\newcommand\norm[1]{\Vert#1\Vert}

\newcommand\abs[1]{\lvert#1\rvert}

\makeatletter
\def\ps@pprintTitle{%
  \let\@oddhead\@empty
  \let\@evenhead\@empty
  \let\@oddfoot\@empty
  \let\@evenfoot\@oddfoot
}

\makeatother
\usepackage[margin=1.00in]{geometry}    
\theoremstyle{definition}
\newtheorem{example}{Example}[section]

\begin{document}

\begin{frontmatter}
  \title{An adaptive dynamically low-dimensional approximation method for multiscale stochastic diffusion equations}

  \author[cuhk]{Eric T. Chung}
  \ead{tschung@math.cuhk.edu.hk}
  \author[cuhk]{Sai-Mang Pun}
  \ead{smpun@math.cuhk.edu.hk}
  \author[hku]{Zhiwen Zhang\corref{cor1}}
  \ead{zhangzw@hku.hk}

  \address[cuhk]{ Department of Mathematics, The Chinese University of Hong Kong,  Shatin, Hong Kong SAR, China.}
  \address[hku]{Department of Mathematics, The University of Hong Kong, Pokfulam Road, Hong Kong SAR, China.}
  \cortext[cor1]{Corresponding author}

\begin{abstract}
In this paper, we propose a dynamically low-dimensional approximation method to solve a class of time-dependent multiscale stochastic diffusion equations. In \cite{ChengHouZhang1:13,ChengHouZhang2:13}, a dynamically bi-orthogonal (DyBO) method was
developed to explore low-dimensional structures of stochastic partial differential equations (SPDEs) and solve them efficiently. However, when the SPDEs have multiscale features in physical space, the original DyBO method becomes expensive. To address this issue, we construct multiscale basis functions within \rev{the framework of generalized multiscale finite element method (GMsFEM)} for dimension reduction in the physical space. To further improve the accuracy, we also perform online procedure to construct online adaptive basis functions. In the stochastic space, we use the generalized polynomial chaos (gPC) basis functions to represent the stochastic part of the solutions.  Numerical results are presented to demonstrate the efficiency of the proposed method in solving time-dependent PDEs with multiscale and random features.
\end{abstract}
\begin{keyword}
Uncertainty quantification (UQ);  dynamically low-dimensional approximation; online adaptive method; stochastic partial differential equations (SPDEs); generalized multiscale finite element method (GMsFEM).
\end{keyword}
\end{frontmatter}
\section{Introduction} \label{sec:introduction}
\noindent
Uncertainty arises in many real-world problems of scientific applications, such as heat propagation through random media or flow driven by stochastic forces. These kind of problems usually have multiple scale features involved in the spatial domain. For example, to simulate flows in heterogeneous porous media, the permeability field is often parameterized by random fields with multiple-scale structures.
 
Stochastic partial differential equations (SPDEs), which contain random variables or stochastic processes, play important roles in modeling complex problems and quantifying the corresponding uncertainties. Considerable amounts of efforts have been devoted to study SPDEs, see \cite{babuska:04,Ghanem:91,WuanHou:06,Zabaras:09,Knio:01,matthies:05,Najm:09,Webster:08,Wan:06,Xiu:09} and references therein. These methods are effective when the dimension of solution space is not huge. However, when SPDEs have multiscale features, the SPDE problems become difficult since it requires tremendous computational resources to resolve the small scales of the SPDE solutions. This motivates us to develop efficient numerical schemes to solve these challenging problems.
 
In this paper, we shall consider the time-dependent SPDEs with multiscale coefficients as follows
\begin{equation}
 \frac{\partial u^{\varepsilon}}{\partial t}(x,t,\omega) = \mathcal{L}^{\varepsilon}u^{\varepsilon} (x,t,\omega), \quad x \in \mathcal{D}, \quad t\in (0,T], \quad \omega \in \Omega, \label{Model_Eq}
\end{equation}
where suitable boundary and initial conditions are imposed, $\mathcal{D}\subset \mathbb{R}^d$ is a bounded spatial domain, $\Omega$ is a sample space, and $\mathcal{L}^{\varepsilon}$ is an elliptic operator that contains multiscale and random coefficient, where the smallest-scale is parameterized by $\varepsilon$. 

The major difficulties in solving \eqref{Model_Eq} come from two parts. In the physical space, we need a mesh fine enough to resolve the small-scale features. In the random space, we need extra degrees of freedom to represent the random features. Moreover, the problem \eqref{Model_Eq} becomes more difficult if the dimension of the random input is high. 
To address these challenges, we shall explore low-dimensional structures hidden in the solution $u^{\varepsilon}(x,t,\omega)$. Specifically, if the solution $u^{\varepsilon}(x,t,\omega)$ is a second-order stochastic process at each time $t>0$, i.e., $u^{\varepsilon}(x,t,\omega) \in L^2(\mathcal{D}\times\Omega)$ for each $t$, one can approximate the solution $u^{\varepsilon}(x,t,\omega)$ by its $m$-term truncated Karhunen-Lo\`{e}ve (KL) expansion \cite{Karhunen1947, Loeve1978}
\begin{equation}
u^{\varepsilon}(x,t,\omega) \approx \bar{u}^{\varepsilon}(x,t) + \sum_{i=1}^m u_i^{\varepsilon}(x,t) Y_i(\omega,t) = \bar{u}^{\varepsilon}(x,t) + \vU (x,t) \vY^T(\omega,t),
\label{KLE}
\end{equation}
where $\vU(x,t) = (u_1^{\varepsilon}(x,t),\cdots,u_m^{\varepsilon}(x,t))$ and $\vY(\omega,t) = (Y_1(\omega,t),\cdots,Y_m(\omega,t))$. The KL expansion gives the compact representation of the solution. However, the direct computation of the KL expansion can be quite expensive since we need to form a covariance kernel and solve a large-scale eigenvalue problem.
 
In \cite{ChengHouZhang1:13,ChengHouZhang2:13}, a dynamically bi-orthogonal (DyBO) method was developed. This new method derives an equivalent system that can faithfully track the KL expansion of the SPDE solution. In other words, the DyBO method gives the evolution equations for $\bar{u}^{\varepsilon}$, $\vU$, and $\vY$. The DyBO method can accurately and efficiently solve many time-dependent SPDEs, such as stochastic Burger's equations and stochastic Navier-Stokes equations, with considerable savings over existing numerical methods. To explore the low-dimensional features of the solutions to the SPDEs, a dynamically orthogonal (DO) method was proposed \cite{sapsis:09}. Later on, the equivalence of DO method and DyBO method has been shown in \cite{choi2014equivalence} and the effectiveness of the DO and DyBO has also been discussed theoretically in \cite{musharbash2015error}. This area is very active and highly demanded due to the latest advances in the UQ research. 

If the SPDEs have multiscale features in the physical space (i.e., the smallest-scale parameter $\varepsilon$ is extremely small), however, the original DyBO method (as well as the DO method) becomes computationally expensive since one needs enormous degrees of freedom to represent the multiscale features in the  physical space. To overcome this difficulty, we shall apply the generalized multiscale finite element methods (GMsFEM) \cite{chung2016adaptive,efendiev2013generalized} to construct multiscale basis functions within each coarse grid block for model reduction in the physical space.

In the GMsFEM, we divide the computation into two stages: the offline  stage and the online  one. In the offline stage, we first compute global snapshot functions within each coarse neighborhood based on the given coarse and fine meshes and construct multiscale basis functions to represent the local heterogeneities. When the snapshot functions are computed, one can construct the multiscale basis functions in each coarse patch by solving some well-designed local spectral problems and identify the crucial multiscale basis functions to form the offline function space. In the online stage, we add more online multiscale basis functions that are constructed using the offline space. These online basis functions are computed adaptively in some selected spatial regions based on the current local residuals and their construction is motivated by the analysis in \cite{chung2015online}. In general, the algorithm  can guarantee that additional online multiscale basis functions will reduce the error rapidly if one chooses a sufficient number of offline basis functions.
We should point out that there are many existing methods in the literature to solve multiscale problems though; see
\cite{Engquist2012,Bourgeat1984, Cruz1995, Dykaar1992,Engquist2003,HanZhangCMS:12,Hou1997,Juanes2005,Peterseim2014,Nolen2008} and references therein. Most of these methods are designed for multiscale problems with deterministic coefficients.
 
In our new method, we first derive the DyBO formulation for the multiscale SPDEs \eqref{Model_Eq}, which consists of deterministic PDEs for $\bar{u}^{\varepsilon}$ and $\vU$ respectively and an ODE system for the stochastic basis $\vY$. For the deterministic PDEs (for $\bar{u}^{\varepsilon}$ and $\vU$) in the formulation, we shall apply the GMsFEM to construct multiscale basis functions and use these multiscale basis functions to represent $\bar{u}^{\varepsilon}$ and $\vU$. It leads to considerable savings over the original DyBO method. For the ODE system, the memory cost is relatively small and we shall apply a suitable ODE solver to compute the numerical solution. The GMsFEM enables us to significantly improve the efficiency of the DyBO method in solving time-dependent PDEs with multiscale and random coefficients.
 
The rest of the paper is organized as follows. In Section \ref{sec:multiscaleSPDEs}, we will introduce the framework of DyBO formulation. The GMsFEM and its online adaptive algorithm will be outlined in Section \ref{sec:multiscale}. The implementation issues of the algorithm and the numerical results will be given in Section \ref{sec:applications}. Finally, some concluding remarks will be drawn in Section \ref{sec:conclusion}.

\section{The DyBO formulation for multiscale time-dependent SPDEs} \label{sec:multiscaleSPDEs}
\noindent
In this paper, we consider a class of parabolic equations with multiscale and random coefficients
\begin{subequations}
\label{eq:MsDyBO_Model}
\begin{align}
\frac{\partial u^{\varepsilon}}{\partial t}
 & =  \mathcal{L}^\varepsilon u^\varepsilon &
~ & x\in \mathcal{D}, ~ t\in (0,T], ~ \omega\in\Omega,  \label{eq:MsDyBO_Model_Eq} \\
u^\varepsilon (x,0,\omega) & =  u_0(x,\omega)  & ~ &x \in \mathcal{D}, ~ \omega \in \Omega, \label{eq:MsDyBO_Model_ic} \\
\mathcal{B}\big( u^\varepsilon(x,t,\omega) \big) & =  h(x, t, \omega) & ~  & x\in \partial \mathcal{D}, ~ \omega \in \Omega. \label{eq:MsDyBO_Model_bc}
\end{align}
\end{subequations}
where $\mathcal{D} \subset \mathbb{R}^{d}$ ($d=2,3$) is a bounded spatial domain, $(\Omega,\mathcal{F}, \mathbb{P}) $ is
a probability space, and suitable boundary and initial conditions are imposed. The differential operator $\mathcal{L}^{\varepsilon}$ is defined as $ \mathcal{L}^\varepsilon u^\varepsilon := \nabla \cdot (a^\varepsilon (x,\omega) )\nabla u^\varepsilon + f(x)$. The multiscale information is described by the parameter $\varepsilon$ and the force $f:\mathbb{R}^d \to \mathbb{R}$ is in $L^2(\mathcal{D})$. Assume that there exist two constants $a_{\max} \gg a_{min}>0$ such that $\mathbb{P}(\omega\in \Omega: a^{\varepsilon}(x,\omega)\in [a_{\min},a_{\max}], \text{a.e. } x \in \mathcal{D}) = 1$. Note that we are interested in the case that the coefficient $a^\varepsilon(x,\omega)$ has high contrast within the domain $\mathcal{D}$, where the mode reduction technique in physical space is necessary to reduce degrees of freedom in representing the solution. 

\subsection{An abstract framework for SPDEs} \label{sec:brief-overview-dybo}
\noindent
To make this paper self-contained, we briefly review the DyBO method \cite{ChengHouZhang1:13,ChengHouZhang2:13}. 
We assume the solution $u^{\varepsilon}(x,t,\omega)$ to \eqref{eq:MsDyBO_Model} satisfies $u^{\varepsilon}(\cdot, t, \cdot) \in L^2(\mathcal{D}\times \Omega)$ for each time $t\in (0,T]$. We consider the $m$-term truncated KL expansion with $m \in \mathbb{N}^+$
\begin{equation}
 \widetilde{u}^{\varepsilon}(x,t,\omega) = \bar{u}^{\varepsilon}(x,t) + \sum_{i=1}^{m} u_i^{\varepsilon}(x,t) Y_i(\omega,t) = \bar{u}^{\varepsilon}(x,t) + \vU(x,t) \vY^T(\omega,t) \approx u^{\varepsilon}(x,t,\omega), \label{eq:KLE:truncated}
\end{equation}
as an approximation to the solution $u^{\varepsilon}(x,t,\omega)$. Here, $\bar{u}^{\varepsilon}(x,t)$ is the mean of the solution,
$$\vU(x,t)=\bkr{u_1^{\varepsilon}(x,t),\cdots,u_m^{\varepsilon}(x,t)} \quad \text{and} \quad \vY(\omega,t)=\bkr{Y_1(\omega,t),\cdots,Y_m(\omega,t)}$$ are the spatial and stochastic modes (with zero mean), respectively. 
We omit the symbol $\varepsilon$ to simplify the notation. Next, by imposing the orthogonal conditions for $\vU$ and $\vY$
$$ \inp{\vU^T}{\vU} = (\inp{u_i}{u_j} \delta_{ij}) \quad \text{and} \quad \EEp{\vY^T \vY} = \mI_{m\times m},$$
we obtain the evolution equations for $\bar{u}^{\varepsilon}$, $\vU$ and $\vY$ as follows
\begin{subequations}
  \label{eq:SPDE:new3}
  \begin{align}
    \diffp{\bar{u}}{t} &= \EEp{\oL \widetilde{u}}, \label{eq:SPDE:new3a} \\
    \diffp{\vU}{t} &= -\vU \mD^T + \EEp{\oLt \widetilde{u} \vY}, \label{eq:SPDE:new3b} \\
    \diffd{\vY}{t} &= -\vY\mC^T + \inp{\oLt \widetilde{u}}{\vU}\mLa_{\vU}^{-1}, \label{eq:SPDE:new3c}
  \end{align}
\end{subequations}
where $\inp{\cdot}{\cdot}$ denotes the inner product in $L^2(\mathcal{D})$, $\mLa_{\vU}=\diag(\inp{\vU^T}{\vU}) \in \RR^{m \times m}$, and $\oLt \widetilde{u} = \oL \widetilde{u} -\EEp{\oL \widetilde{u}}$. Define two operators $\oQ: \RR^{k \times k} \goto \RR^{k \times k}$ and $\oQt: \RR^{k \times k} \goto \RR^{k \times k}$ as follows
\begin{equation*}
  \oQ(\mM)  := \frac{1}{2}\bkr{\mM-\mM^T}\quad \text{and} \quad \oQt(\mM) := \oQ(\mM) + \diag(\mM),
\end{equation*}
where $\mM \in \RR^{k \times k}$ is a square matrix and $\diag(\mM)$ is a diagonal matrix whose diagonal entries are equal to those of $\mM$. Then, the matrices $\mC,\mD \in \RR^{m\times m}$ in \eqref{eq:SPDE:new3} can be solved uniquely from the following linear system
\begin{subequations}
  \label{eq:CDSystem}
  \begin{align}
    \mC - \mLa_{\vU}^{-1} \oQt\bkr{\mLa_{\vU}\mC} &= 0, \label{eq:CDSystem:C}\\
    \mD-\oQ\bkr{\mD} &= 0, \label{eq:CDSystem:D} \\
    \mD^T+\mC &= G_{*}(\bar{u}^{\varepsilon},\vU,\vY), \label{eq:CDSystem:CD}
  \end{align}
\end{subequations}
where the matrix is given by $G_{*}(\bar{u}^{\varepsilon},\vU,\vY)=\mLa_{\vU}^{-1}\inp{\vU^T}{\EEp{\oLt \widetilde{u}^{\varepsilon} \vY}} \in \RR^{m \times m}$.
 
In order to represent the stochastic modes $Y_i(\omega,t)$ in \eqref{eq:SPDE:new3c}, one can choose several different approaches including ensemble representations in sampling methods and spectral representations. 
In this paper, we use gPC basis functions to represent the stochastic modes $Y_i(\omega,t)$. 
Given two positive integers $r$ and $p$, we define $\sJ_r^p := \{ \minda: \minda = (\alpha_1,\cdots, \alpha_r ), \alpha_i \in \mathbb{N}_+, \abs{\alpha} = \sum_{i=1}^r \alpha_i \leq p \} \backslash \{0 \}$.
Let $\{ H_i(\xi) \}_{i=1}^\infty$ denote as an one-dimensional family of $\rho$-orthogonal polynomial, i.e., 
$$\int_{-\infty}^\infty H_i(\xi) H_j(\xi) \rho(\xi) \ d\xi = \delta_{ij}.$$
If we write $\hmpn_{\minda}(\vxi) = \prod_{i=1}^r H_{\alpha_i} (\xi_i)$ for $\minda \in \sJ_r^p$ and $\vxi = (\xi_i)_{i=1}^r$, then the Cameron-Martin theorem \cite{Cameron1947} implies the stochastic modes $Y_i(\omega,t)$ in \eqref{eq:KLE:truncated} can be approximated by 
\begin{equation}
  \label{eq:YgPC:component}
  Y_i(\omega,t) \approx \sum_{\minda \in \sJ_r^p} \hmpn_{\minda}(\vxi(\omega))A_{\minda i}(t) = \hmpn\bkr{\vxi} \mA_i(t) \quad i=1,2,\cdots,m.
\end{equation}
Here, $\hmpn\bkr{\vxi}=\bkr{\hmpn_{\minda}\bkr{\vxi}}_{\minda \in \sJ_r^p} \in \mathbb{R}^{1\times N_p}$, 
$\mA_i(t) = (\mA_{\minda i}(t))_{\minda \in \sJ_r^p} \in \mathbb{R}^{N_p \times 1}$ and $N_p := \abs{\sJ_r^p}$. 
\rev{We remark that for each $i = 1,2,\cdots,m$, the coefficients $\{ \mA_{\minda i}(t) \}_{\minda \in \sJ_r^p}$ represent the projection coefficients of the stochastic mode $Y_i(\omega,t)$ on the gPC basis functions
$\hmpn_{\minda}(\vxi)$, $\minda \in \sJ_r^p$. Moreover, $\{ \mA_{\minda i}(t) \}_{\minda \in \sJ_r^p}$ change with respect to time. }
One may write
\begin{equation}
  \label{eq:YgPC:matrix}
  \vY(\omega, t) = \hmpn\bkr{\vxi(\omega)}\mA (t)
\end{equation}
where $\mA(t) = (\mA_1(t), \cdots, \mA_m(t)) \in \RR^{N_p \times m}$. The KL expansion \eqref{eq:KLE:truncated} now reads
\begin{equation*}
\widetilde{u} \approx \bar{u}+\vU\mA^T\hmpn^T.
\end{equation*}
We can derive equations for $\bar{u}$, $\vU$ and $\mA$, instead of $\bar{u}$, $\vU$ and $\vY$. 
\rev{Here and in the following, we have suppressed the variables $x$, $t$, and $\omega$ for notation simplicity.} In other words, the stochastic modes $\vY$ are identified with a matrix $\mA \in \RR^{N_p \times m}$, which leads to the DyBO-gPC formulation of  SPDE \eqref{eq:MsDyBO_Model}
\begin{subequations}
  \label{eq:DyBO:gPC}
  \begin{align}
    \diffp{\bar{u}}{t} &= \EEp{\oL \widetilde{u}}, \label{eq:DyBO:gPC:a} \\
    \diffp{\vU}{t} &= -\vU \mD^T + \EEp{\oLt \widetilde{u} \hmpn}\mA, \label{eq:DyBO:gPC:b} \\
      \diffd{\mA}{t} &= -\mA \mC^T + \inp{\EEp{\hmpn^T\oLt \widetilde{u}}}{\vU}\mLa_{\vU}^{-1}, \label{eq:DyBO:gPC:c}
  \end{align}
\end{subequations}
where $\mC(t)$ and $\mD(t)$ can be solved from the linear system \eqref{eq:CDSystem} with
\begin{equation}
  \label{eq:DyBO:gPC:Gs}
  G_{*}(\bar{u},\vU,\vY)=\mLa_{\vU}^{-1}\inp{\vU^T}{\EEp{\oLt \widetilde{u} \vY}} = \mLa_{\vU}^{-1}\inp{\vU^T}{\EEp{\oLt \widetilde{u} \hmpn}}\mA.
\end{equation}
By solving the system \eqref{eq:DyBO:gPC}, we have an approximate solution to \eqref{eq:MsDyBO_Model}
\begin{equation*}
  u^{\text{DyBO-gPC}}=\bar{u}+\vU\mA^T\hmpn^T.
\end{equation*}
The condition $\EEp{\vY^T \vY}$ implies that the columns $(\mA_i)_{i=1}^m$ are orthonormal, i.e., $\mA^T\mA= \mI_{m\times m}$. Note that $\mA\mA^T \in \RR^{N_p \times N_p}$ in general is not an identity matrix as $m \ll N_p$ if the SPDE solution has a low-dimensional structure.
\subsection{The DyBO formulation for the model problem} \label{sec:analysis}
\noindent
In this section, we shall derive the DyBO formulation for the model problem \eqref{eq:MsDyBO_Model}. Recall that 
the definition of the differential operator is $\mathcal{L}u = \nabla \cdot (a(x,\omega) \nabla u) + f(x)$ and we have omitted $\varepsilon$ for notation simplification.  We assume that the coefficient $a(x,\omega)$ is of the form
$a(x,\omega)=\bar{a}(x)+\tilde{a}(x,\omega)$,
where $\bar{a}(x) = \EEp{a(x,\omega)}$ and $\tilde{a}(x,\omega)$ is the fluctuation, which can be parametrized as follows
\begin{equation*} 
\tilde{a}(x,\omega) = \sum_{i=1}^r a_i(x) \xi_i(\omega) , \quad i = 1,2\cdots,r.
\end{equation*}
Here, $r\geq1$ is a positive integer and $\{ \xi_i(\omega)\}_{i=1}^r$ are independent identically distributed random variables assumed to be mean-zero.
By substituting the expression of $\mathcal{L}u$ into \eqref{eq:DyBO:gPC}, we obtain the DyBO-gPC formulation for
the model problem \eqref{eq:MsDyBO_Model} (see \ref{cha:dybo-formulation-MsSPDE} for the details of the derivation)
\begin{subequations}
\label{eq:DyBO-gPC-MsFEM}
\begin{align}
  \diffp{\bar{u}}{t}
  & =  \nabla\cdot(\bar{a}\nabla \bar{u}) + \nabla\cdot(\EEp{\tilde{a}\nabla \vU \mA^{T} \hmpn^T})+ f, \label{eq:DyBO-gPC-MsFEM-ubar}\\
  \diffp{\vU}{t}
  &= -\vU\mD^T +  \nabla\cdot(\EEp{\tilde{a} \nabla \bar{u}\hmpn})\mA
+   \nabla\cdot(\bar{a} \nabla \vU) + \nabla\cdot( \EEp{ \tilde{a} \nabla \vU \mA^{T} \hmpn^T \hmpn})\mA, \label{eq:DyBO-gPC-MsFEM-U}\\
  \diffd{\mA}{t}
  &= -\mA \mC^T + \inp{ \nabla\cdot(\EEp{\hmpn^{T}\tilde{a}\nabla\bar{u}})+
  \nabla\cdot(\bar{a} \nabla \mA\vU^{T})  + \nabla\cdot( \EEp{ \tilde{a} \hmpn^T \hmpn \mA \nabla \vU^{T} })  }{\vU}\mLa_{\vU}^{-1}, \label{eq:DyBO-gPC-MsFEM-A}
\end{align}
\end{subequations}
where matrices $\mC$ and $\mD$ can be solved from \eqref{eq:CDSystem} with  $G_{*} = \mLa_{\vU}^{-1}\inp{\vU^T}{\EEp{\oLt u \hmpn}}\mA$.

\begin{remark}
When the force function of  the model problem \eqref{eq:MsDyBO_Model} contains randomness, i.e., $f(x,\omega)$, 
one can derive the DyBO formulation accordingly without any difficulty. 
\end{remark}
\begin{remark} \label{rmk:GeneralizationToSystems}
The boundary conditions and initial conditions for each physical component, and the initial condition for each stochastic component can be obtained by projection of the initial and boundary conditions of $u(x,t,\omega)$ on the corresponding components. 
\end{remark}

\begin{remark}
As the system evolves, the norm of the mode $u_i$ (denoted as $\lambda_i$) in the KL expansion may change and some of them may get closer to each other. In this case, if the matrices $\mC$ and $\mD$ are still solved from \eqref{eq:CDSystem}, numerical errors will pollute the results. One may freeze $\vU$ or $\vY$ temporarily for a short time and continue to evolve the system. At the end of this short period, the solution is recast into the bi-orthogonal form via the KL expansion. See \cite[Section 4.2]{ChengHouZhang1:13} for more details.
\end{remark}

\section{Multiscale model reduction using the GMsFEM} \label{sec:multiscale}
\subsection{Motivations}\label{sec:motivations}
\noindent
Since the DyBO formulation \eqref{eq:DyBO-gPC-MsFEM} involves multiscale features in the physical space, one may consider an efficient solver to solve the problem in order to reduce the computational cost. As such, we shall apply the GMsFEM to discretize $\bar{u}$ and $\vU$. Note that, Eq.\eqref{eq:DyBO-gPC-MsFEM-ubar} and each component of Eq. \eqref{eq:DyBO-gPC-MsFEM-U} have the following deterministic time-dependent PDE form
\begin{subequations}
\label{eq:GMsFEM}
    \begin{align}
        \frac{\partial w}{\partial t} & = \nabla \cdot ( \bar{a} \nabla w) + \mathcal{G}, \label{eq:GMsFEM:1}\\
        w|_{t=0}  & = w_0. \label{eq:GMsFEM:2}
    \end{align}
for some functions  $\mathcal{G}$. For example, in \eqref{eq:DyBO-gPC-MsFEM-ubar} we have $w=\bar{u}$ and $\mathcal{G} = \nabla \cdot (\EEp{\tilde{a} \nabla \vU \mA^{T} \hmpn^{T}}) + f$.
\end{subequations}
 
In order to discretize the equation \eqref{eq:GMsFEM} in time, we apply the implicit Euler scheme with time step $\Delta t>0$ and obtain the discretization for each time $t_n = n \Delta t$, $n=1,2,\cdots, N$ ($T = N\Delta t$)
\begin{equation*}
    \frac{w^n - w^{n-1}}{\Delta t} = \nabla \cdot(\bar{a} \nabla w^n) + \mathcal{G}
\end{equation*}
where $w^n = w(t_n)$ and the above equation is equivalent to the following
\begin{equation}
    -\nabla \cdot (\bar{a}\nabla w^n) + cw^n =  \tilde{\mathcal{G}}, \label{eq:GMsFEM:diff:1}
\end{equation}
where $c = 1/\Delta t$ and $\tilde{\mathcal{G}} = cw^{n-1} + \mathcal{G}$. Hence, for each fixed $t_n>0$, we use the GMsFEM to solve the second order elliptic PDE \eqref{eq:GMsFEM:diff:1} with multiscale coefficient $\bar{a}$.
\subsection{The GMsFEM and the multiscale basis functions}
\noindent
Next, we present the framework of the GMsFEM for solving \eqref{eq:GMsFEM:diff:1}. We first introduce the notion of fine and coarse grids that we shall use in the method. Let $\mathcal{T}^H$ be a conforming partition of the spatial domain $\mathcal{D}$ with mesh size $H>0$. We refer to this partition as the coarse grid. Subordinate to $\mathcal{T}^H$, we define a fine grid partition denoted by $\mathcal{T}^h$, with mesh size $0<h \ll H$, by refining each coarse element in $\mathcal{T}^H$ into a connected union of fine elements. 
Assume the above refinement is performed such that $\mathcal{T}^h$ is a conforming partition of $\mathcal{D}$.
Denote the interior nodes of $\mathcal{T}^H$ as $x_i$, $i = 1,2,\cdots, N_{in}$,
where $N_{in}$ is the number of interior nodes. The coarse elements of $\mathcal{T}^H$ are denoted as $K_j$, $j=1,2.,\cdots,N_{e}$, where $N_{e}$ is the number of the coarse elements. Define the coarse neighborhood of the node $x_i$ by $D_i:=\bigcup \{ K_j \in \mathcal{T}^H: x_i \in \overline{K_j} \}$. 
 
Once the coarse and fine grids are given, one may construct the multiscale basis functions for approximating the solution of \eqref{eq:GMsFEM:diff:1}. To obtain the multiscale basis functions, we first define the snapshot space.
For each neighborhood $D_i$, define $J_h(D_i)$ as the set of fine nodes of $\mathcal{T}^h$ lying on $\partial D_i$ and denote its cardinality as $L_i \in \mathbb{N}^+$. For each fine-grid node $x_j \in J_h(D_i)$, 
define a fine-grid function $\delta_j^h$ on $J_h(D_i)$ as $\delta_j^h(x_k) = \delta_{jk}$. Next, for $j = 1,\cdots, L_i$, define the snapshot function $\psi_j^{(i)}$ in coarse neighborhood $D_i$ as the solution to the following system
\begin{eqnarray}
    -\nabla \cdot(\bar{a} \nabla \psi_j^{(i)})  & = 0 &\quad \text{in } D_i, \\
    \psi_j^{(i)} & =  \delta_j^h &\quad \text{on } \partial D_i.
\end{eqnarray}
The local snapshot space $V_{snap}^{(i)}$ corresponding to the coarse neighborhood $D_i$ is defined as follows $V_{snap}^{(i)} := \text{snap} \{ \psi_j^{(i)}: j = 1,\cdots, L_i \}$ and the snapshot space reads $V_{snap} := \bigoplus_{i=1}^{N_{in}} V_{snap}^{(i)}$.
 
The snapshot space defined above is usually of large dimension. Therefore, a dimension reduction is performed on $V_{snap}$ to archive a smaller space $V_{\text{off}}$, which contains the multiscale basis functions for simulation. This reduction is achieved by performing a spectral decomposition on each local snapshot space $V_{snap}^{(i)}$. The analysis in \cite{Yalchin2011} motivates the following construction. For each $i = 1,\cdots,N_{in}$, the spectral problem is to find $(\phi_j^{(i)},\lambda_j^{(i)}) \in V_{snap}^{(i)} \times \mathbb{R}$ such that
\begin{eqnarray}
    \int_{D_i} \bar{a} \nabla \phi_j^{(i)} \cdot \nabla v  = \lambda_j^{(i)} \int_{D_i} \hat{a} \phi_j^{(i)} v \quad \forall v\in V_{snap}^{(i)}, \quad j = 1,\cdots, L_i, \label{eqn:spectral}
\end{eqnarray}
where $\hat{a} := \bar{a} \sum_{i=1}^{N_{in}} H^2 \abs{\nabla \chi_i}^2$ and $\{\chi_i\}_{i=1}^{N_{in}}$ is a set of partition of unity satisfying the following system 
    \begin{eqnarray*}
        -\nabla \cdot (\bar{a} \nabla \chi_i)  &= 0  & \quad \text{in } K\subset D_i, \label{eqn:std_basis_1}\\
        \chi_i &= p_i &  \quad \text{on each }  \partial K ~ \text{with } K \subset D_i \label{eqn:std_basis_2}, \\
        \chi_i &= 0 & \quad \text{on } \partial D_i. \label{eqn:std_basis_3}
    \end{eqnarray*}
Assume that the eigenvalues obtained from \eqref{eqn:spectral} are arranged in ascending order and we may use the first $0<l_i \leq L_i$ (with $l_i \in \mathbb{N}^+$) eigenfunctions (related to the smallest $l_i$ eigenvalues) to form the local multiscale space $V_{\text{off}}^{(i)} := \text{span} \{ \chi_i \phi_j^{(i)} : j = 1,\cdots, l_i\}$. The multiscale space $V_{\text{off}}$ is the direct sum of the local multiscale spaces, namely $V_{\text{off}} := \bigoplus_{i=1}^{N_{in}} V_{\text{off}}^{(i)}$.
 
Once the multiscale space $V_{\text{off}}$ is constructed, we can find the GMsFEM solution $u_{\text{off}}^n \in V_{\text{off}}$ at time $t=t_n$, $n = 1,\cdots, N$, by solving the following equation
\begin{equation}
    \mathcal{A}(u_{\text{off}}^n,v) + c\inp{u_{\text{off}}^n}{v} = \inp{c u_{\text{off}}^{n-1} + \mathcal{G}}{v} \quad \forall v \in V_{\text{off}}, \label{eqn:GMsFEM:var}
\end{equation}
where $\mathcal{A}(u,v) := \int_{\mathcal{D}} \bar{a} \nabla u \cdot \nabla v $.
\begin{remark}
The above derivation of $V_{\text{off}}$ is based on the (mean) coefficient $\bar{a}$ and the multiscale basis functions in $V_{\text{off}}$ are suitable for approximating $\bar{u}$. In this paper, we assume that the fluctuation of the coefficient is a small perturbation to the mean. Therefore, the multiscale space $V_{\text{off}}$ can also efficiently approximate $\vU$.
\end{remark}
\subsection{Online adaptive algorithm}
\noindent
In order to achieve a rapid convergence in the GMsFEM, one may add some online basis functions to enrich the multiscale space $V_{\text{off}}$ based on local residuals. In this subsection, we briefly outline the online adaptive algorithm for GMsFEM.
 
Let $u_{\text{off}}^n \in V_{\text{off}}$ be the numerical solution obtained in \eqref{eqn:GMsFEM:var} at time $t=t_n$. Given a coarse neighborhood $D_i$, we define $V_i := H_0^1(D_i) \cap V_{snap}$ equipped with the norm $\norm{v}_{V_i}^2 := \int_{D_i} \bar{a}(x) \abs{\nabla v}^2$. We also define the local residual operator $\mathcal{R}_i^n : V_i \to \mathbb{R}$ by
\begin{equation}
    \mathcal{R}_i^n(v;u_{\text{off}}^n) := \int_{D_i} \big(cu_{\text{f}}^{n-1} + \mathcal{G}\big)v - \int_{D_i} \big(\bar{a} \nabla u_{\text{off}}^n \cdot \nabla v + cu_{\text{off}}^n v \big) \quad \forall v \in V_i,\label{eq:residualonline}
\end{equation}
where $u_{\text{f}}^{n-1}$ is the fine-scale solution at time $t=t_{n-1}$. The operator norm of $\mathcal{R}_i^n$, denoted by $\norm{\mathcal{R}_i^n}_{V_i^*}$ gives a measure of the quantity of the residual. The online basis functions are computed during the time-marching process for a given fixed time $t = t_n$, contrary to offline basis functions that are pre-computed (defined in the Section 4.2).
 
Suppose that one needs to add an online basis $\phi$ into the space $V_i$. 
The analysis in \cite{chung2015online} suggests that the required online basis $\phi \in V_i$ is the solution to the following equation
\begin{equation}
\mathcal{A}(\phi,v) = \mathcal{R}_i^n(v;u_{\text{off}}^{n,\tau}) \quad \forall v\in V_i. \label{eqn:online}
\end{equation}
We refer to $\tau \in \mathbb{N}$ as the level of enrichment and denote $u_{\text{off}}^{n,\tau}$ as the solution of \eqref{eqn:GMsFEM:var} in $V_{\text{off}}^{n,\tau}$. Remark that $V_{\text{off}}^{n,0} := V_{\text{off}}$ for time level $n \in \mathbb{N}$. Let $\mathcal{I} \subset \{1,2,\cdots,N_{in}\}$ be the index set over some non-overlapping coarse neighborhoods. For each $i\in \mathcal{I}$, we obtain an online basis $\phi_i \in V_i$ by solving \eqref{eqn:online} and define $V_{\text{off}}^{n,\tau+1} = V_{\text{off}}^{n,\tau} \oplus \text{span}\{ \phi_i: i\in \mathcal{I}\}$. After that, solve \eqref{eqn:GMsFEM:var} in $V_{\text{off}}^{n,\tau+1}$ to get $u_{\text{off}}^{n,\tau+1}$. 
Consequently, following the arguments in \cite{chung2015online}, we have at time $t=t_n$,
\begin{equation}
\norm{u_{\text{f}}^n - u_{\text{off}}^{n,\tau+1}}_V^2 \leq \bigg( 1- \frac{\Lambda_{\text{min}}^{(\mathcal{I})}}{C_{\text{err}}}  \frac{\sum_{i\in \mathcal{I}} \norm{\mathcal{R}_i^n}_{V_i^*}(\lambda_{l_i+1}^{(i)})^{-1} }{\sum_{i=1}^N \norm{\mathcal{R}_i^n}_{V_i^*}(\lambda_{l_i+1}^{(i)})^{-1}} \bigg) \norm{u_{\text{f}}^n-u_{\text{off}}^{n,\tau}}_V^2, \label{ieq:online}
\end{equation}
where $C_{\text{err}}$ is a uniform constant and $\Lambda_{\text{min}}^{(\mathcal{I})} = \min_{i\in \mathcal{I}} \lambda_{l_i +1}^{(i)}$. Here, the norm is defined by $\norm{\cdot}_V := \sqrt{\mathcal{A}(\cdot,\cdot)}$.
Inequality \eqref{ieq:online} shows that we can obtain a better accuracy by adding more online basis functions at each time $t=t_n$ and the rate of convergence depends on the constant $C_{\text{err}}$ and $\Lambda_{\text{min}}^{(\mathcal{I})}$.

\subsection{The implementation of our new algorithm}\label{sec:ImplementationGMsFEM}
\noindent
We summarize the computational scheme for the problem in this section. Recall that the multiscale coefficient is $a(x,\omega) = \bar{a}(x) + \tilde{a}(x,\omega)$. The mean has high contrast in nature and the fluctuation part is $\tilde{a}(x,\omega) = \sum_{i=1}^r a_i \xi_i = a_i \xi_i$, where the Einstein notation is used. We rewrite the DyBO formulation for \eqref{eq:DyBO:gPC} as follows 
\begin{align}
  \diffp{\bar{u}}{t}
  & =  \nabla\cdot(\bar{a}\nabla \bar{u}) + \nabla\cdot(\EEp{a_i \xi_i \nabla \vU \mA^{T} \hmpn^T})+ f, \label{eqn:DyBO-gPC-ubar}\\
  \diffp{\vU}{t}
  &= -\vU\mD^T + \nabla\cdot(\EEp{a_i \xi_i \nabla \bar{u}\hmpn})\mA
+   \nabla\cdot(\bar{a} \nabla \vU) + \nabla\cdot( \EEp{ a_i \xi_i \nabla \vU \mA^{T} \hmpn^T \hmpn})\mA, \label{eqn:DyBO-gPC-U}\\
  \diffd{\mA}{t}
  &= -\mA \mC^T + \inp{ \nabla\cdot(\EEp{\hmpn^{T} a_i \xi_i \nabla\bar{u}})+
  \nabla\cdot(\bar{a} \nabla \mA\vU^{T})  + \nabla\cdot( \EEp{ a_i \xi_i \hmpn^T \hmpn \mA \nabla \vU^{T} })  }{\vU}\mLa_{\vU}^{-1}. \label{eqn:DyBO-gPC-A}
\end{align}
We assume that homogeneous boundary condition is imposed. Hence, the solutions $\bar{u}$ and $\vU = (u_1,\cdots,u_m)$ will vanish on $\partial \mathcal{D}$. If the model problem \eqref{eq:MsDyBO_Model} has inhomogeneous boundary condition, the boundary conditions for $\bar{u}$ and $\vU$ can be obtained 
by taking expectation of $u$ on the corresponding components. The initial conditions for $\bar{u}$, $\vU$ and $\mA$ depend on the initial condition of $u$, which will be discussed in Section \ref{sec:applications}. For the details of implementation, see \ref{app:implement}.

\section{Numerical experiments} \label{sec:applications}
\noindent
In this section, we present some numerical examples to demonstrate the efficiency of our proposed method. The computational domain is $\mathcal{D}=(0,1)^2 \subset \mathbb{R}^2$ and $T = 1$. First, we divide the domain $\mathcal{D}$ into several equal square units with mesh size $H>0$ and refer to it as the coarse mesh $\mathcal{T}^H$. Next, we divide each coarse element into several equal square blocks with mesh size $h>0$ and refer to it as the fine mesh $\mathcal{T}^h$. Then, we discretize $\bar{u}$ and $\vU$ in the DyBO formulation 
by using the GMsFEM to reduce degrees of freedom in representing the multiscale solutions. Thus, with the multiscale basis functions, we can represent $\bar{u}$ and $\vU$  on the coarse mesh. \rev{In all examples, the number of initial local basis functions is $L_i = 4$}. 
 
The multiscale coefficient is assumed to be $a(x,\omega) = \bar{a}(x) + \sum_{i=1}^r a_i(x)\xi_i(\omega)$, where $a_i(x)$ is a small (or multiscale) perturbation and $\{\xi_i(\omega)\}_{i=1}^r$ is a set of i.i.d. uniform-distributed random variables over $[-1,1]$. Moreover,  we assume that there exist two constants $a_{\max} \gg a_{\min}>0$ such that $\mathbb{P}(\omega\in \Omega: a(x, \omega)\in [a_{\min},a_{\max}], ~\text{a.e. } x \in \mathcal{D}) = 1$.

The initial condition of the solution is assumed to have the form of $m$-terms truncated KL expansion
\begin{equation}
\tilde{u}(x,0,\omega) = \bar{u}(x,0) + \sum_{i=1}^m u_i(x,0) Y_i(\omega,0).\label{eq:MsDyBO_Initial_KLEform}
\end{equation}
The stochastic basis $Y_i(\omega,t)$ can be expanded as
$ Y_i(\omega,t) = \sum_{j=1}^{N_p} H_j(\omega) A_{ji}(t)$ for each $i = 1,\cdots, m$. Here, $\{ H_j(\omega)\}_{j=1}^{N_p}$ is a set of tensor products of orthogonal polynomials in $\mathbb{R}$, $N_p = \frac{(p+r)!}{p!r!}-1$, and $p$ is the maximum degree of polynomial. Denote $\mA(t) = (A_{ji}(t))_{N_p \times m}$. The initial condition of the matrix $\mA(t)|_{t=0} = \big(A_{ji}(0)\big)_{N_p \times m}$ should satisfy $ \EEp{\hmpn \mA} = 0$ and $\mA^T(t) \mA (t) |_{t=0} = \mI_{m \times m}$.
 
For each function to be approximated (e.g. $\bar{u}$, $u_i$ or the variance function $\text{var}(u) := \sum_{i=1}^m u_i^2$), we define the following quantity at $t= t_n$ to measure the numerical error 
$$e_2^n = \frac{\norm{u_{\text{f}}^n - u_{\text{approx}}^n}_{L^2(\mathcal{D})}}{\norm{u_{\text{f}}^n}_{L^2(\mathcal{D})}},
$$
where $u_{\text{f}}^n$ is the reference solution and $u_{\text{approx}}^n$ is the approximation obtained by the proposed method. In the remaining part of this paper, we refer to this quantity as $L^2$-error.

\begin{example}\label{exp:1}
\noindent
We set the mesh size to be $H = \sqrt{2}/10$ and $h = \sqrt{2}/100$. The time step is $\Delta t = 10^{-3}$. The multiscale fluctuation is parameterized by three independent random variables ($r=3$) and the number of terms in the KL expansion is $m=4$. 
 Next, we set the coefficients $a_i$ ($i=1,2,3$) to be 
\begin{align*}
  a_1(x_1,x_2) & = 0.04 \times \frac{2+P_1 \sin(\frac{2\pi (x_1-x_2)}{\varepsilon_1})}{2-P_1 \cos(\frac{2\pi (x_1-x_2)}{\varepsilon_1})}, & P_1 = 1.6 \quad \text{and} \quad \varepsilon_1 = 1/8, \\
a_2 (x_1,x_2) & = 0.08 \times \frac{2+P_2 \cos(\frac{2\pi x_1}{\varepsilon_2})}{2-P_2 \sin(\frac{2\pi x_2}{\varepsilon_2})}, & P_2 = 1.5 \quad \text{and} \quad \varepsilon_2 = 1/7, \\
a_3(x_1,x_2) & = 0.16 \times \frac{2+P_3 \sin(\frac{2\pi (x_1-0.5)}{\varepsilon_3})}{2-P_3 \cos(\frac{2\pi (x_2-0.5)}{\varepsilon_3})}, & P_3 = 1.4 \quad \text{and} \quad \varepsilon_3 = 1/6.
\end{align*}
The mean $\bar{a}$ of the multiscale coefficient is of high-contrast (see Figure \ref{fig:exp3_abar}). The source function is chosen to be $ f \equiv 1$. The initial conditions for the mean of the solution and the physical modes are given as follows
\begin{eqnarray*}
 \bar{u}(x_1,x_2,t)|_{t=0} & = & 32\big(1-\cos(2\pi x_1)\big) \big(1-\cos(2\pi x_2) \big), \\
 u_{1}(x_1,x_2,t)|_{t=0} & = & 24 \big(1-\cos(2\pi x_1)\big) \big(1-\cos(2\pi x_2) \big), \\
 u_{2}(x_1,x_2,t)|_{t=0} & = & 16 \big(1-\cos(4\pi x_1)\big) \big(1-\cos(4\pi x_2) \big), \\
 u_{3}(x_1,x_2,t)|_{t=0} & = & 8 \big(1-\cos(6\pi x_1)\big) \big(1-\cos(6\pi x_2) \big), \\
 u_{4}(x_1,x_2,t)|_{t=0} & = & 4 \big(1-\cos(8\pi x_1)\big) \big(1-\cos(8\pi x_2) \big).
\end{eqnarray*}

\begin{figure}[htbp!]
\mbox{
\begin{subfigure}{0.49\textwidth}
\centering
\includegraphics[width=0.9\linewidth, height=5.2cm]{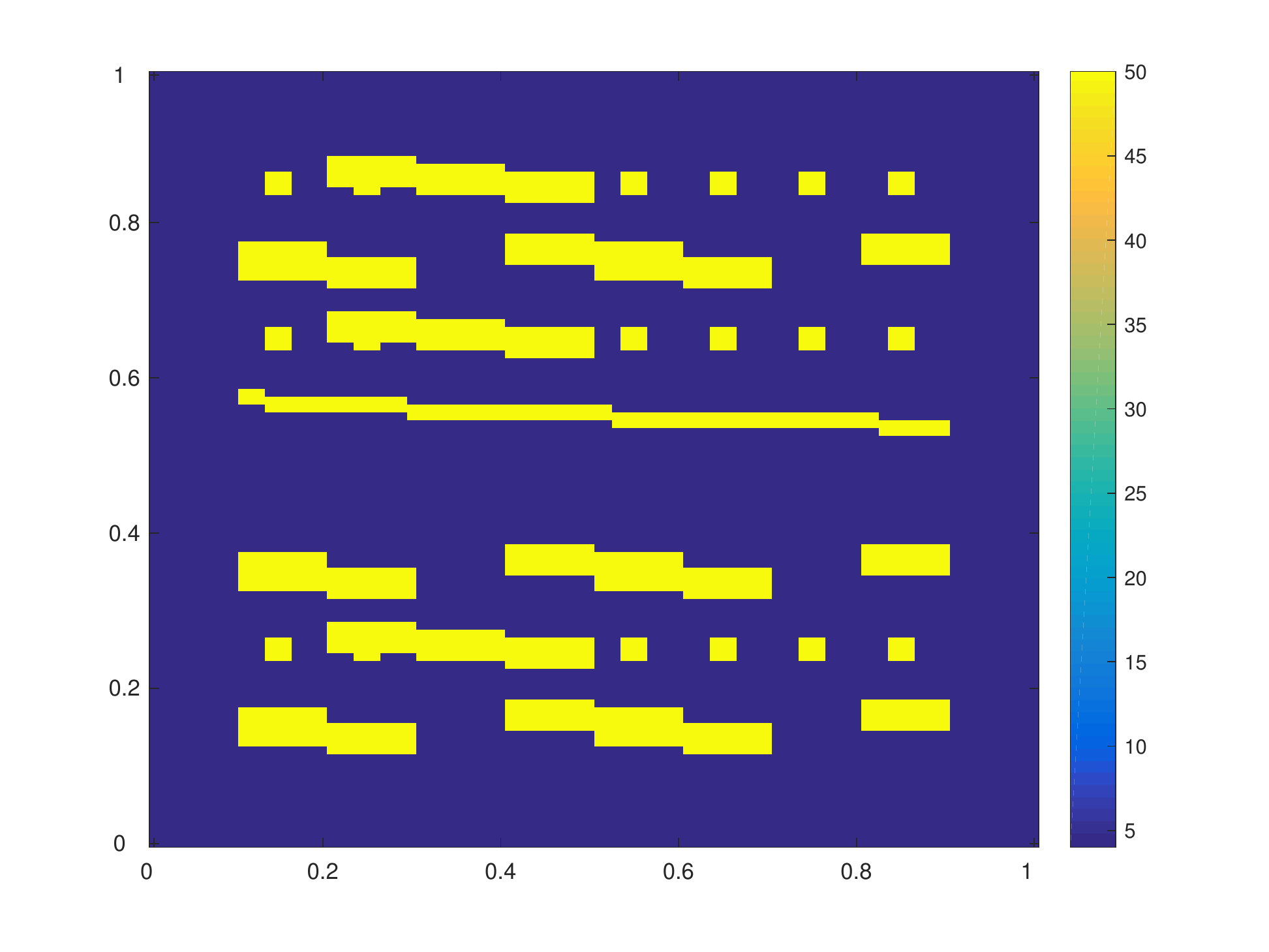} 
\caption{$\bar{a}$ in the Example \ref{exp:1}. (Max: \rev{1000}, Min: 4)}
\label{fig:exp3_abar}
\end{subfigure}

\begin{subfigure}{0.49\textwidth}
\centering
\includegraphics[width=0.9\linewidth, height=5.2cm]{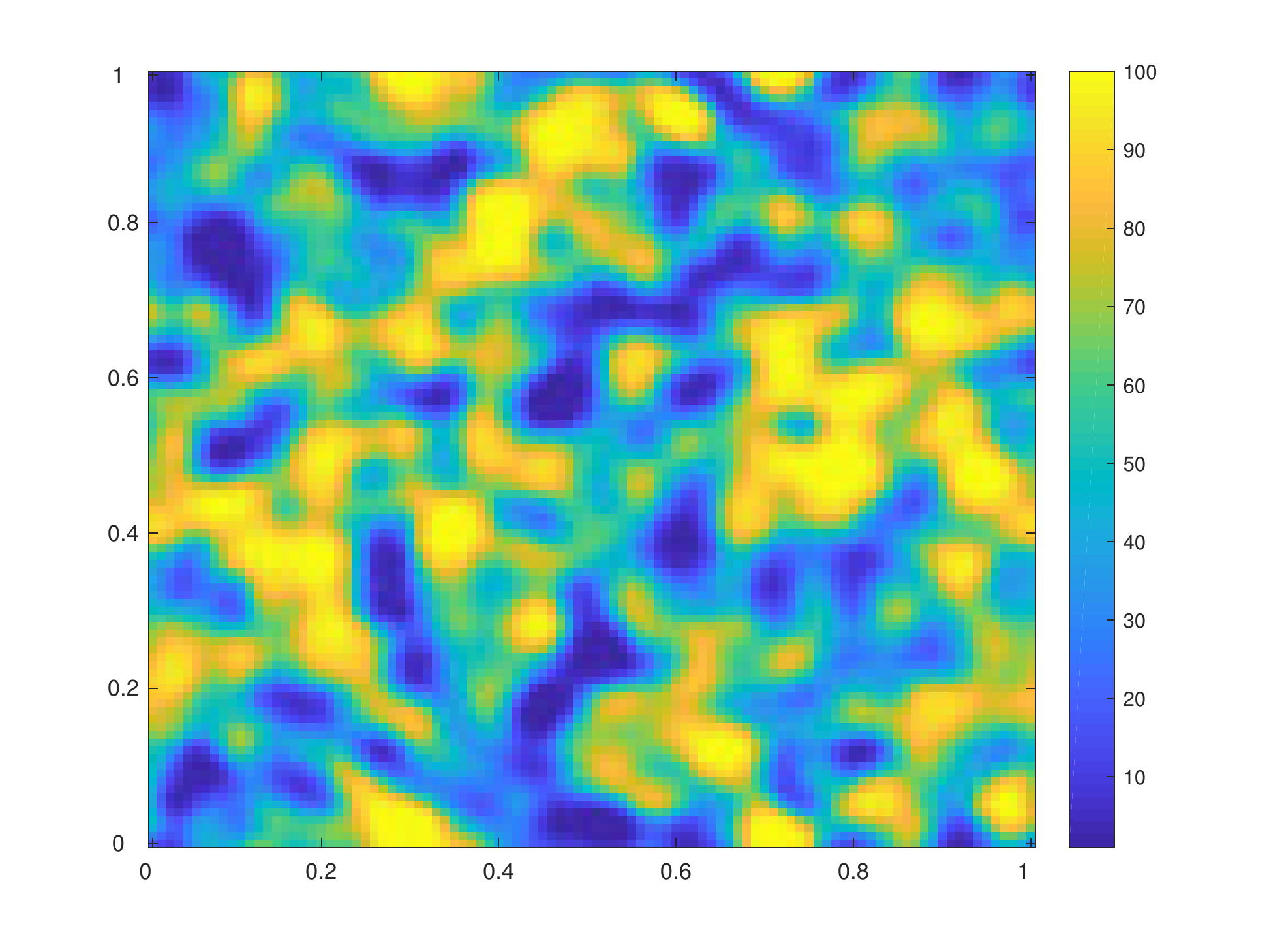}
\caption{$\bar{a}$ in the Example \ref{exp:2}. (Max: 100, Min: 1)}
\label{fig:exp_spe10}
\end{subfigure}
}
\caption{The mean component of the permeability.}
\label{fig:mean_exp}
\end{figure}

The history of the $L^2$-error is recorded in Table \ref{tab:exp3_2}. 
One can find that at the specific time level the $L^2$-errors of the quantities to be approximated are relative small (less than $1\%$) when the online procedure is terminated. It shows that the proposed method can approximate the stochastic multiscale diffusion problem with certain accuracy.
We remark that due to the linearity of the diffusion problem and its DyBO formulation, one may easily extend this algorithm to the case with more modes in the KL expansion. In addition, one can adopt the adaptive approach proposed
in \cite{ChengHouZhang2:13} to dynamically change the number of the modes in the DyBO formulation during the numerical simulation.


\begin{table}[h!]
\centering
\begin{tabular}{c|c|ccccc}
\hline
\hline
function & online status & $t=0.1$   & $t=0.2$   & $t=0.4$   & $t=0.8$   & $t=1.0$   \\
\hline
\multirow{2}{*}{$\bar{u}$}       & S & 3.7904\% & 4.0551\% & 4.0404\% & 4.0278\% & 4.0378\% \\
                                 & E & 0.3883\%  & 0.4076\%  & 0.4059\%  & 0.4042\%  & 0.4062\%  \\
\hline
\multirow{2}{*}{$u_1$}           & S & 3.7449\% & 5.0027\% & 4.3672\% & 4.2385\% & 4.5720\% \\
                                 & E & 0.4247\%  & 0.4000\%  & 0.3489\%  & 0.3964\%  & 0.3641\%  \\
\hline
\multirow{2}{*}{$u_2$}           & S & 4.9319\% & 6.7831\% & 7.5986\% & 5.0301\% & 5.2988\% \\
                                 & E & 0.4284\%  & 0.4496\%  & 0.5077\%  & 0.3662\%  & 0.3708\%  \\
\hline
\multirow{2}{*}{$u_3$}           & S & 8.2879\% & 14.5263\% & 5.7078\% & 14.5705\% & 8.4522\% \\
                                 & E & 0.5011\% &    0.5505\% &    0.4391\% &    0.5294\% &    0.4251\% \\
\hline
\multirow{2}{*}{$u_4$}           & S & 14.9988 \% &   12.6509\% &   11.0844\% &   22.7257\% &   20.3689\%  \\
                                 & E & 0.6415\% &    0.4641\% &    0.4586\% &    0.6168\% &    0.6943   \% \\                                                               
\hline
\multirow{2}{*}{$\text{var}(u)$} & S & 6.3057\% & 7.6541\% & 6.7792\% & 6.6286\% & 6.8721\% \\
                                 & E & 0.6795\%  & 0.6729\%  & 0.5870\%  & 0.6181\%  & 0.6258\% \\
\hline
\hline
\end{tabular}
\caption{$L^2$-error for each functions in Example \ref{exp:1}. (S: start, E: end)}
\label{tab:exp3_2}
\end{table}
\end{example}

\begin{example}\label{exp:2}
We keep $H$ and $h$ the same as in Example \ref{exp:1}. The time step is still $\Delta t = 10^{-3}$. The mean permeability field $\bar{a}$ in this example is chosen from the SPE10 data set \cite{aarnes2005mixed} and the data is moderately related to the real physical applications (See Figure \ref{fig:exp_spe10}). The fluctuation part is parameterized by four independent random variables ($r=4$) and the coefficients is set as $a_i(x_1,x_2)=0.02 \times \frac{2+P_i\sin(\frac{2\pi (x_1-x_2)}{\varepsilon_i})}{2-P_i\cos(\frac{2\pi(x_1+x_2)}{\varepsilon_i})}$, $i=1,...,4$, where 
$[P_1,P_2,P_3,P_4]=[1.4,1.5,1.6,1.7]$ and $[\varepsilon_1,\varepsilon_2,\varepsilon_3,\varepsilon_4]=[\frac{1}{9},\frac{1}{8},\frac{1}{7},\frac{1}{6}]$. 
The number of modes in the KL expansion is $m=3$ and the source function is $f \equiv 1$. The initial conditions for the mean and the physical modes are given as follows
\begin{eqnarray*}
 \bar{u}(x_1,x_2,t)|_{t=0} & = & 4\big(1-\cos(2\pi x_1)\big) \big(1-\cos(2\pi x_2) \big), \\
 u_{1}(x_1,x_2,t)|_{t=0} & = & 16 \big(1-\cos(4\pi x_1)\big) \big(1-\cos(4\pi x_2) \big), \\
 u_{2}(x_1,x_2,t)|_{t=0} & = & 4 \big(1-\cos(6\pi x_1)\big) \big(1-\cos(6\pi x_2) \big), \\
 u_{3}(x_1,x_2,t)|_{t=0} & = & 2 \big(1-\cos(8\pi x_1)\big) \big(1-\cos(8\pi x_2) \big).
\end{eqnarray*}
\begin{figure}[h!]
\mbox{
\begin{subfigure}{0.5\textwidth}
\centering
\includegraphics[width=0.82\linewidth, height = 4.9cm]{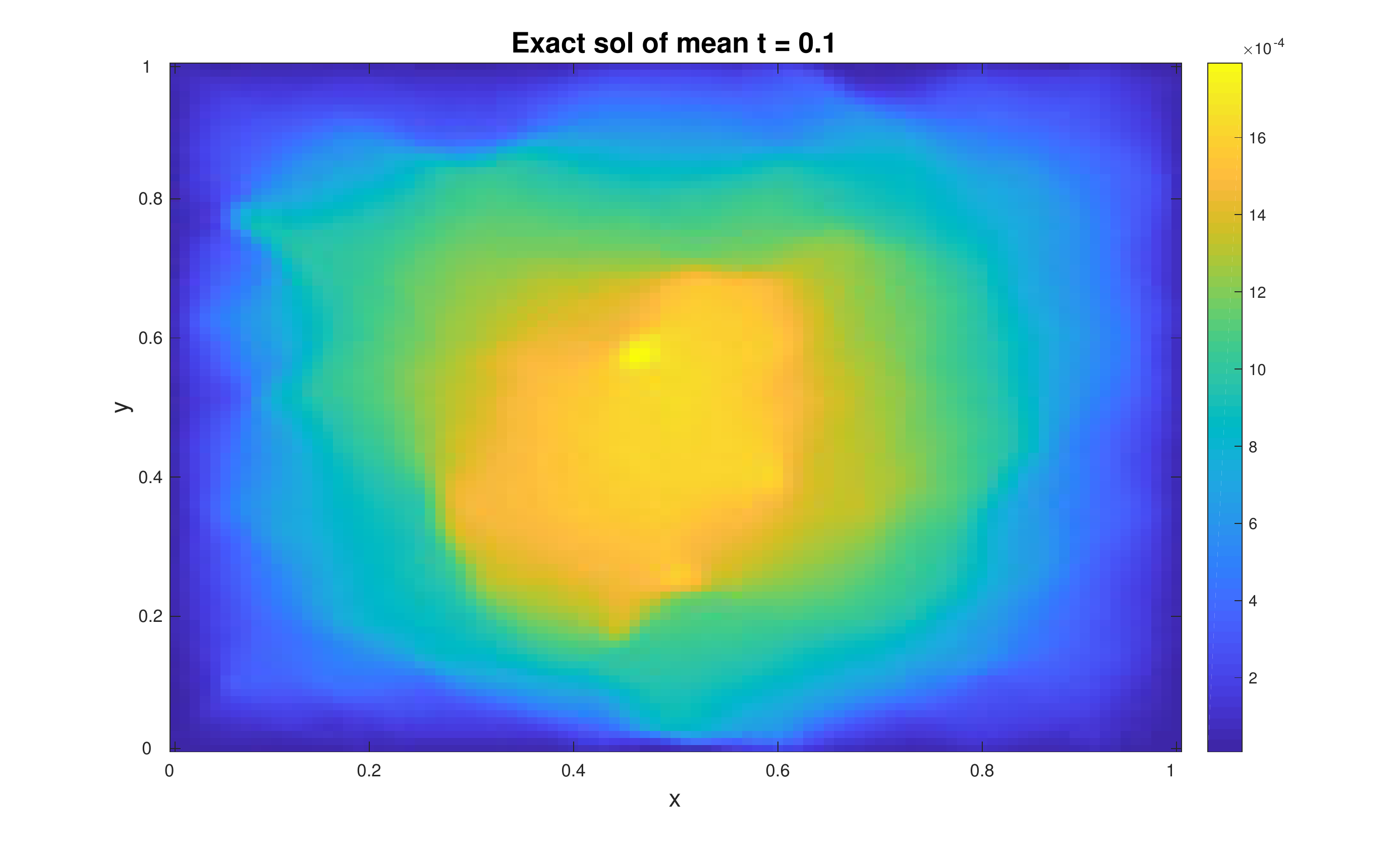}
\caption{Fine-scale solution of the mean $\bar{u}$.}
\label{fig:fs_mean}
\end{subfigure}

\begin{subfigure}{0.5\textwidth}
\centering
\includegraphics[width=0.82\linewidth, height = 4.9cm]{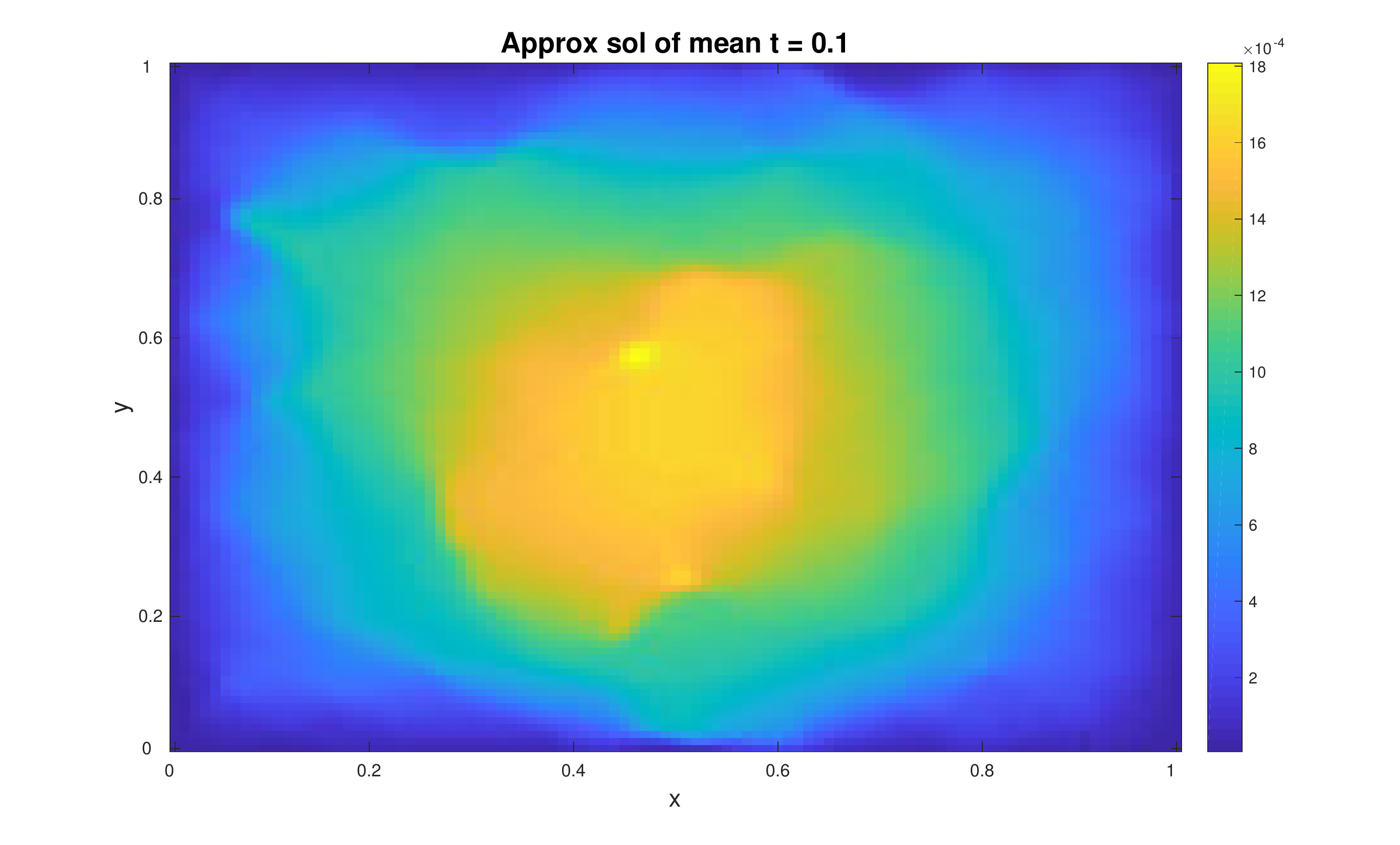}
\caption{Multiscale approximation of the mean $\bar{u}$.}
\label{fig:ms_mean}
\end{subfigure}
}
\mbox{
\begin{subfigure}{0.5\textwidth}
\centering
\includegraphics[width=0.82\linewidth, height = 4.9cm]{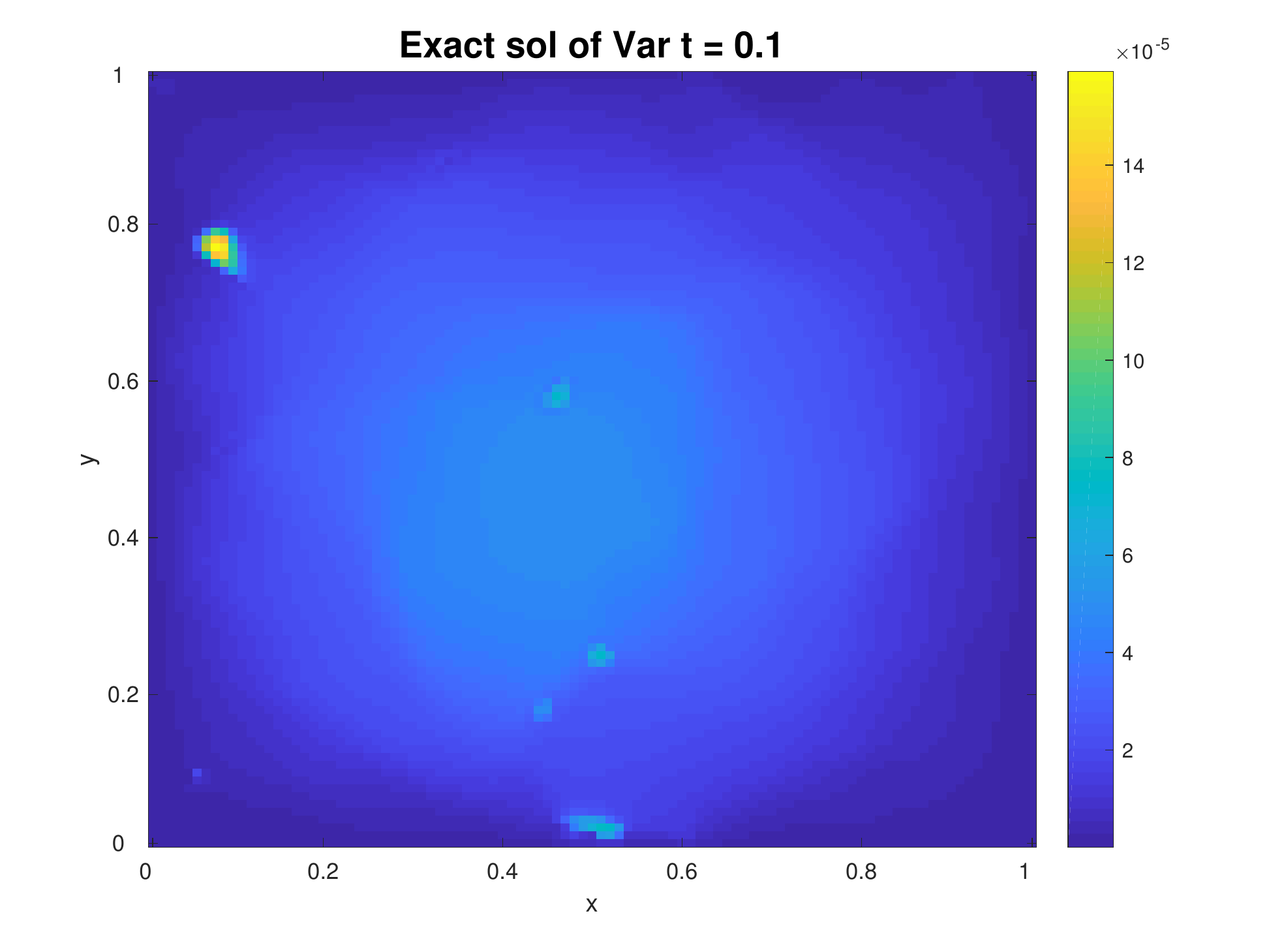}
\caption{Fine-scale solution of the variance.}
\label{fig:fs_mode1}
\end{subfigure}

\begin{subfigure}{0.5\textwidth}
\centering
\includegraphics[width=0.82\linewidth, height = 4.9cm]{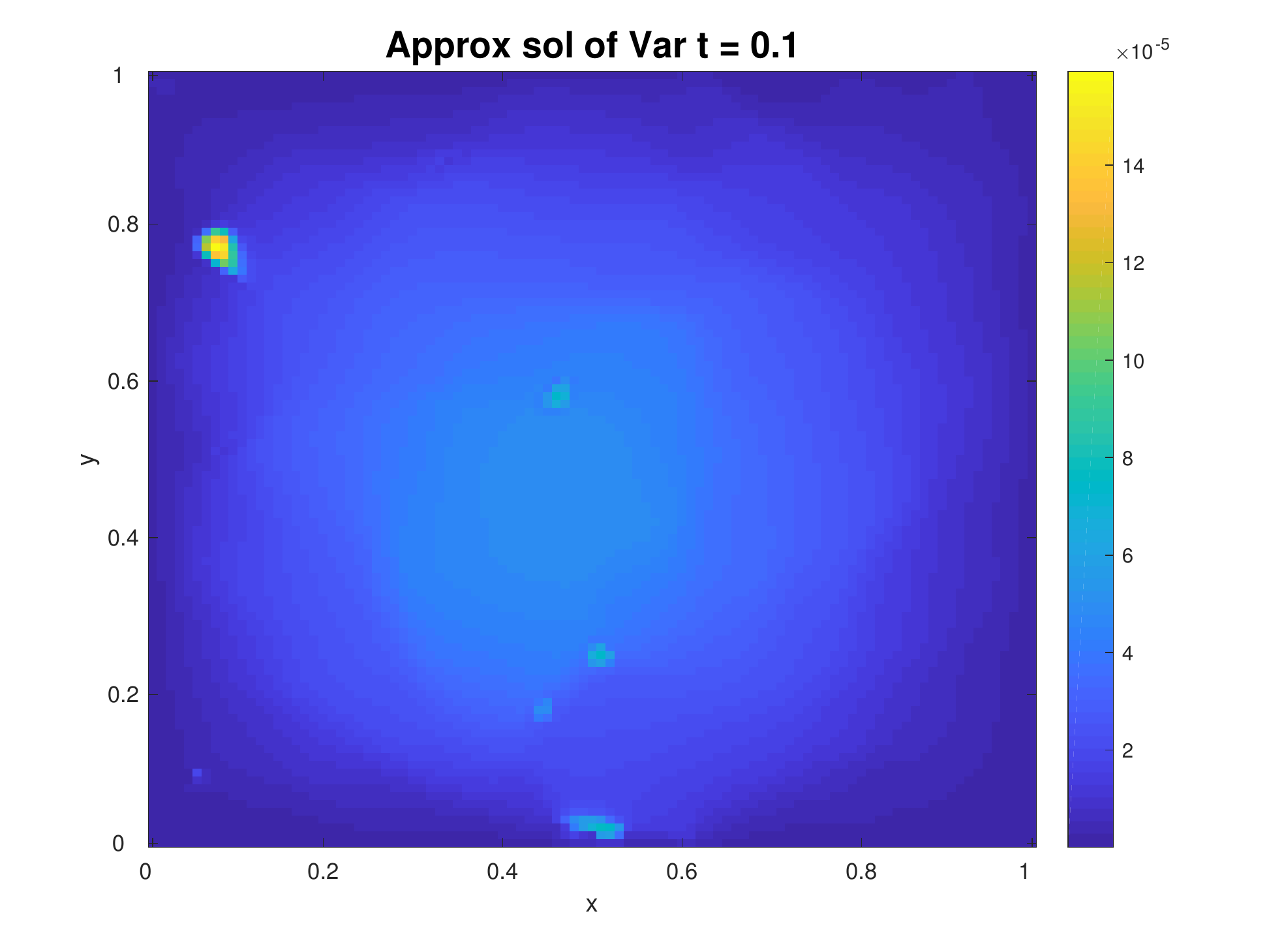}
\caption{Multiscale approximation of the variance.}
\label{fig:ms_mode1}
\end{subfigure}
}
\caption{Solution profiles at $t=0.1$ in Example \ref{exp:2}.}
\label{fig:sol_profiles}
\end{figure}
One may notice that this problem has multiscale features driven by the mean field $\bar{a}$ and some small random perturbations. The solution profiles of the mean and the variance at $t=0.1$ are plotted in Figure \ref{fig:sol_profiles}. One can see that our method archives a certain level of accuracy when the problem has both multiscale and random features. 

In both the numerical experiments, a few times of online enrichments are required at each time level. Meanwhile, the $L^2$-error between the multiscale solution and the fine-scale solution is nearly less than $2\%$ when the online procedure is terminated.
 
We remark that the contrast value in SPE10 model used in Example \ref{exp:2} is already scaled down by 100 times. The difficulty of these kinds of stochastic multiscale problem is that when the contrast value is high (e.g. $\max_{x \in \mathcal{D}} \big(a(x)\big)\approx 10^4$ or larger), the usual computational schemes for UQ problems fail to obtain a good approximation, even though the random perturbation is small. We shall develop a more robust method to compute stochastic multiscale problems with higher contrast value in our subsequent research.

\end{example}

\begin{example}\label{exp:3}
\rev{In this example, we compare the efficiency between the proposed method and the fine-scale method in terms of the CPU time. First, we set the mesh size to be $H = \sqrt{2}/10$ and $h = \sqrt{2}/400$. The time step is $\Delta t = 1/80$ and the final time is $T=1$. We keep the setting of random variables and initial conditions the same as in Example \ref{exp:1}. Table \ref{tab:exp_CPU_1} records the $L^2$-error at certain time level using online adaptivity. We remark that when the multiscale space contains sufficiently many initial local basis functions, only $1$-$2$ times of iterations are required to achieve such certain accuracy. Furthermore, one may achieve moderate computational savings with the proposed multiscale solver. From the data in Table \ref{tab:exp_CPU_2}, one may observe that the proposed DyBO-GMsFEM solver outperform the fine-scale solver with $20$ times speed-up in terms of the CPU times.}  

\begin{table}[h!]
\centering
\begin{tabular}{c|ccccc}
\hline
\hline
function & $t=1/8$   & $t=1/4$   & $t=1/2$   & $t=3/4$   & $t=1$   \\
\hline
$\bar{u}$ & 0.2286\% &    0.2317\% &    0.2621\% &    0.1922\% &    0.3450\%  \\

$u_1$ & 0.2416\% &    0.2758\% &    0.2659\% &    0.2497\% &    0.2478\%  \\

$u_2$ & 0.2774\% &    0.3218\% &    0.3403\% &    0.2975\% &    0.2791\%  \\

$u_3$ & 0.5053\% &    0.3110\% &    0.3701\% &    0.2835\% &    0.4562\%  \\

$u_4$ & 0.7377\% &    0.2931\% &    0.4040\% &    0.3612\% &    0.4286\%  \\

$\text{var}(u)$ & 0.3911\% &    0.4730\% &    0.4276\% &    0.4219\% &    0.4131\% \\
\hline
\hline
\end{tabular}
\caption{$L^2$-error for each functions in Example \ref{exp:3}; online process terminated.}
\label{tab:exp_CPU_1}
\end{table}

\begin{table}[h!]
\centering
\begin{tabular}{l|cc}
\hline
\hline
 & \multicolumn{2}{c}{CPU times (s)} \\
\hline
function & Fine-scale solver &    Proposed solver \\
\hline
Mean $\bar{u}$ & 44.1012 &    2.0763	\\
Modes $u_i$ & 170.9146 &    9.4423 	\\
Total & $215.0158$ & $11.5186$ 		\\
\hline
\hline
\end{tabular}
\caption{CPU times of the fine-scale solver and the proposed solver. ($H = 1/10$, $h = 1/400$)}
\label{tab:exp_CPU_2}
\end{table}
\end{example}

\section{Conclusion} \label{sec:conclusion}
\noindent
In this paper, we proposed a new framework combining the DyBO formulation and the online adaptive GMsFEM to solve time-dependent PDEs with multiscale and random features. For a given multiscale PDE with random input, one can derive its corresponding DyBO formulation under the assumption that the solution has a low-dimensional structure in the sense of Karhunen-Lo\`{e}ve expansion. The DyBO method enables one to faithfully track the KL expansion of the SPDE solution. 
 
For the mean of the solution and physical modes of the solution in the truncated KL expansion, they are deterministic and dependent on time, which were solved using the GMsFEM with implicit Euler scheme. Moreover, at each time level, the online construction was applied in order to reduce the $L^2$-error rapidly. For the stochastic modes of the solution in the truncated KL expansion, we projected them onto a set of polynomial chaos to obtain an ODE system, which could be solved using some existing solvers. Thanks to the approximation property of the multiscale basis functions obtained using the GMsFEM, the degrees of freedom of our new method is relatively small compared with the original DyBO method. Therefore, our new method provides considerable computational savings over the original DyBO method.

We presented several numerical examples for 2D stochastic parabolic PDEs with multiscale coefficients to demonstrate the accuracy and efficiency of our proposed method. 
\rev{One may obtain significant saving in computation with the proposed multiscale solver without losing the accuracy of approximations}. 
We point out that the stochastic multiscale problem is still very challenging when the contrast value of the coefficient is very large, which will be our subsequent research work.
 

\section*{Acknowledgements}
\noindent
The research of E. Chung is supported by Hong Kong RGC General Research Fund (Project 14304217). The research of Z. Zhang is supported by the Hong Kong RGC grants (Projects 27300616, 17300817, and 17300318), National Natural Science Foundation of China (Project 11601457), and Seed Funding Programme for Basic Research (HKU). S. Pun would like to thank the Isaac Newton Institute for Mathematical Sciences for support and hospitality during the programme {\it Uncertainty quantification for complex systems: theory and methodologies} when work on this paper was undertaken. This workshop was supported by EPSRC Grant Numbers EP/K032208/1 and EP/R014604/1.
\appendix
\section{Derivations of the DyBO Formulation for the multiscale SPDE} \label{cha:dybo-formulation-MsSPDE}
\noindent
In this appendix, we provide the details of the derivations of the DyBO-gPC formulation of multiscale SPDE \eqref{eq:MsDyBO_Model}. Substituting the KL expansion of $u$ (see Eq.~\eqref{KLE}) into Eq.~\eqref{eq:MsDyBO_Model}, we get
\begin{align*}
  \oL u & = \nabla\cdot((\bar{a}+\tilde{a})(\nabla \bar{u} + \nabla \vU \mA^{T} \hmpn^T  )) + f \\
  & = \nabla\cdot(\bar{a} \nabla \bar{u}) + \nabla\cdot(\tilde{a} \nabla \bar{u}) +
   \nabla\cdot(\bar{a} \nabla \vU \mA^{T} \hmpn^T) + \nabla\cdot(\tilde{a} \nabla \vU \mA^{T} \hmpn^T)
   + f.
\end{align*}
Taking expectations on both sides yields
\begin{align*}
  \EEp{\oL u}& = \nabla\cdot(\bar{a}\nabla \bar{u}) + \nabla\cdot(\EEp{\tilde{a}\nabla \vU \mA^{T} \hmpn^T})+ f,
\end{align*}
where we have used the facts that $\EEp{\tilde{a}}=0$ and  $\EEp{\hmpn}=\textbf{0}$. Then, we obtain
\begin{align*}
\oLt u & = \oL u - \EEp{\oL u}  \\
       & =  \nabla\cdot(\tilde{a} \nabla \bar{u})
+ \nabla\cdot(\bar{a} \nabla \vU \mA^{T} \hmpn^T) + \nabla\cdot(\tilde{a} \nabla \vU \mA^{T} \hmpn^T)
- \nabla\cdot(\EEp{\tilde{a}\nabla \vU \mA^{T} \hmpn^T})
\end{align*}
In addition, we compute some related terms as follows
\begin{align*}
  \EEp{\oLt u \hmpn} & = \nabla\cdot(\EEp{\tilde{a} \nabla \bar{u}\hmpn})
+   \nabla\cdot(\bar{a} \nabla \vU \mA^{T}) + \nabla\cdot( \EEp{ \tilde{a} \nabla \vU \mA^{T} \hmpn^T \hmpn})
\end{align*}
and
\begin{align*}
  \inp{\vU^T}{\EEp{\oLt u \hmpn}}_{m \times N_p}
  &=  \inp{\vU^T}{ \nabla\cdot(\EEp{\tilde{a}\nabla\bar{u}\hmpn})+\nabla\cdot(\bar{a} \nabla \vU \mA^{T})  + \nabla\cdot( \EEp{ \tilde{a} \nabla \vU \mA^{T} \hmpn^T \hmpn}) } \\
\end{align*}
From Eq.~\eqref{eq:DyBO:gPC}, we obtain the DyBO-gPC formulation for the multiscale SPDE \eqref{eq:MsDyBO_Model}
\begin{align*}
  \diffp{\bar{u}}{t}
  & =  \nabla\cdot(\bar{a}\nabla \bar{u}) + \nabla\cdot(\EEp{\tilde{a}\nabla \vU \mA^{T} \hmpn^T})+ f, \\
  \diffp{\vU}{t}
  &= -\vU\mD^T +  \nabla\cdot(\EEp{\tilde{a} \nabla \bar{u}\hmpn})\mA
+   \nabla\cdot(\bar{a} \nabla \vU) + \nabla\cdot( \EEp{ \tilde{a} \nabla \vU \mA^{T} \hmpn^T \hmpn})\mA, \\
  \diffd{\mA}{t}
  &= -\mA \mC^T + \inp{ \nabla\cdot(\EEp{\hmpn^{T}\tilde{a}\nabla\bar{u}})+
  \nabla\cdot(\bar{a} \nabla \mA\vU^{T})  + \nabla\cdot( \EEp{ \tilde{a} \hmpn^T \hmpn \mA \nabla \vU^{T} })  }{\vU}\mLa_{\vU}^{-1} ,
\end{align*}
where matrices $\mC$ and $\mD$ can be solved from \eqref{eq:CDSystem} with  $G_{*}$
\begin{align*}
 G_{*} = \mLa_{\vU}^{-1}\inp{\vU^T}{\EEp{\oLt u \hmpn}}\mA.
\end{align*}
and we have used that $\mA^{T}\mA= \mI_{m\times m}$.

\section{The implementation of DyBO-GMsFEM}\label{app:implement}
\noindent
In this section, we present the details of the implementation of our complete algorithm. 
We denote $V_{\text{off}} = \text{span} \{ \eta_i: i = 1,\cdots, N_d \}$ and the row vector $\mathcal{C} = \mathcal{C}(x) = \big (\eta_1(x),\cdots, \eta_{N_d}(x)\big)$, where $N_d = \text{dim}(V_{\text{off}})$. For each time $t = t_n$, we seek the approximations for the functions $\bar{u}$ and $\vU$ using the multiscale basis functions and the following representations hold
$$ \bar{u}(x,t) = \mathcal{C}(x) \hat{u}_0(t), \quad \hat{u}_0(t) \in \mathbb{R}^{N_d},$$
$$ \vU(x,t) = \mathcal{C}(x) \hat{U}_m(t), \quad \hat{U}_m(t) := (\hat{u}_1(t),\cdots,\hat{u}_m(t)) \in \mathbb{R}^{N_d\times m}.$$
Then, the variational form of \eqref{eqn:DyBO-gPC-ubar} becomes
\begin{equation}
\oM \frac{d\hat{u}_0}{dt} = -\oS_0 \hat{u}_0 - \oS_i \hat{U}_m \mA^T \EEp{\xi_i \hmpn^T} + \hat{f}, \label{eq:DyBO-GMsFEM-ubar}
\end{equation}
where $$ \mathcal{M} = ( \inp{\eta_j}{\eta_k} )\in \mathbb{R}^{N_d \times N_d}, \quad \mathcal{S}_0 = (\inp{\bar{a} \eta_j}{\eta_k} )\in \mathbb{R}^{N_d \times N_d},$$
$$ \mathcal{S}_i = ( \inp{a_i \eta_j}{\eta_k})\in \mathbb{R}^{N_d \times N_d}, ~ ~ i=1,...,r, \quad \hat{f} = (\inp{f}{\eta_1} \cdots \inp{f}{\eta_{N_d}})^T \in \mathbb{R}^{N_d}.$$
Similarly, the variational form of \eqref{eqn:DyBO-gPC-U} becomes
\begin{equation}
\oM \frac{d\hat{U}_m}{dt} = - \oM\hat{U}_m \mD^T - \oS_i\hat{u}_0  \EEp{\xi_i \hmpn}\mA - \oS_0 \hat{U}_m - \oS_i \hat{U}_m \mA^T \EEp{\xi_i \hmpn^T \hmpn }\mA,   \label{eq:DyBO-GMsFEM-U}
\end{equation}
where the Einstein notation is used. Next, we apply the implicit Euler method to approximate the time derivatives in \eqref{eq:DyBO-GMsFEM-ubar} and \eqref{eq:DyBO-GMsFEM-U}. Combining with the variational forms, we obtain the following algebraic equations at each fixed time $t = t_n = n\Delta t$, $n = 1,\cdots, N$
\begin{align}
    \oS_0 \hat{u}_0^n + c \oM \hat{u}_0^n & = \mathcal{G}_1^{n-1}, \label{eq:DyBO-GMsFEM-ubar-Euler} \\
    \oS_0 \hat{U}_i^n + c \oM \hat{U}_i^n & = \mathcal{G}_2^{n-1},  \quad i=1,...,m, \label{eq:DyBO-GMsFEM-U-Euler}
\end{align}
where $ c = 1/\Delta t$ and the right hand sides $\mathcal{G}_1$ and $\mathcal{G}_2$ are defined as follows
\begin{align*}
    \mathcal{G}_1^{n-1} & = c \oM \hat{u}_0^{n-1} - \oS_i \hat{U}_m^{n-1} \mA_{n-1}^T \EEp{\xi_i \hmpn^T} + \hat{f}, \\
    \mathcal{G}_2^{n-1} & = c \oM \hat{U}_m^{n-1} - \oM \hat{U}_m^{n-1} \mD_{n-1}^T - \oS_i \hat{u}_0^{n-1}\EEp{\xi_i \hmpn}\mA_{n-1} -\oS_i \hat{U}_m^{n-1}\mA_{n-1}^T \EEp{\xi_i \hmpn^T \hmpn}\mA_{n-1},
\end{align*}
where $\mA_{n-1} = \mA(t_{n-1})$, $\hat{U}_m^n = \hat{U}_m(t_n)$, $\hat{u}_0^n = \hat{u}_0(t_n)$, and $\mD_n = \mD(t_n)$.  Using integration by part one simplifies the ODE system for \rev{$\mA(t) = (\mA_1(t), \cdots , \mA_m(t)) \in \mathbb{R}^{N_p \times m}$} as follows
\begin{align}
\frac{d\mA}{dt} = -\mA \mC^T - \big( \EEp{\xi_i \hmpn^T}\hat{u}_0^T \oS_i \hat{U}_m + \mA \hat{U}_m^T \oS_0 \hat{U}_m + \EEp{\xi_i \hmpn^T \hmpn} \mA \hat{U}_m^T \oS_i \hat{U}_m \big) \mLa_{\vU}^{-1}.
\label{eq:DyBO-GMsFEM-A}
\end{align}
\rev{Here, $\mA_i(t) = (\mA_{\minda i} (t) )_{\minda \in \sJ_r^p} \in \mathbb{R}^{N_p \times 1}$, $i = 1,\cdots, m$ represent the stochastic components of the solution, which change with respect to time.} 
Then, we use implicit Euler scheme to approximate the time derivative and get
\begin{align}
\mA_{n} = \mA_{n-1} - \Delta t \big( \mA_{n-1} \mC_{n-1}^T + \mathcal{G}_3^{n-1}\big),
\label{eq:DyBO-GMsFEM-A-Euler}
\end{align}
where $\mC_{n-1} = \mC(t_{n-1})$ and
$$ \mathcal{G}_3^{n-1} =\Big( \EEp{\xi_i \hmpn^T}(\hat{u}_0^{n-1})^T \oS_i \hat{U}_m^{n-1} + \mA_{n-1}(\hat{U}_m^{n-1})^T \oS_0 \hat{U}_m^{n-1} + \EEp{\xi_i \hmpn^T \hmpn} \mA_{n-1} (\hat{U}_m^{n-1})^T \oS_i \hat{U}_m^{n-1} \Big) \mLa_{\vU}^{-1}.$$

\noindent
Overall, we solve the following discrete system to obtain $\hat{u}_0^n$, $\hat{U}_m^n$, and $\mA_{n}$ at each time $t=t_n$, $n=1,\cdots,N$,
\begin{align}
    \oS_0 \hat{u}_0^n + c \oM \hat{u}_0^n & = \mathcal{G}_1^{n-1}, \label{eqn:final_1st}\\
    \oS_0 \hat{U}_m^n + c \oM \hat{U}_m^n & = \mathcal{G}_2^{n-1}, \label{eqn:final_2nd}\\
    \mA_{n} & = \mA_{n-1} - \Delta t \big( \mA_{n-1} \mC_{n-1}^T + \mathcal{G}_3^{n-1}\big), \label{eqn:final_3rd}
\end{align}
where the matrices $\mC_{n-1}$ and $\mD_{n-1}$ in \eqref{eqn:final_1st}-\eqref{eqn:final_3rd}
can be computed using the system \eqref{eq:CDSystem} with $G_{*}(\bar{u},\vU,\vY) = -\mLa_{\vU}^{-1} \big( \hat{U}_m^T \oS_i^T \hat{u}_0 \EEp{\xi_i \hmpn} + \hat{U}_m^T \oS \hat{U}_m \mA^T + \hat{U}_m^T \oS_i^T \hat{U}_m \mA^T \EEp{\xi_i \hmpn^T \hmpn}\big) \mA $.

To improve the accuracy of the spatial approximation, one may possibly perform the online adaptive enrichment at each time level, adjusting the dimension of the multiscale space. See \cite{chung2015online} for more details of the online basis construction using GMsFEM.


\bibliographystyle{plain}
\bibliography{mypaperref}

\end{document}